\documentclass[11pt]{article}
\usepackage{latexsym,amssymb,amsmath}
\begin{document}

\newcommand{\rd}{\mbox{Rad}}
\newcommand{\kn}{\mbox{ker}}
\newcommand{\psp}{\vspace{0.4cm}}
\newcommand{\pse}{\vspace{0.2cm}}
\newcommand{\ptl}{\partial}
\newcommand{\dlt}{\delta}
\newcommand{\Dlt}{\Delta}
\newcommand{\si}{\sigma}
\newcommand{\cta}{\theta}
\newcommand{\al}{\alpha}
\newcommand{\be}{\beta}
\newcommand{\G}{\Gamma}
\newcommand{\g}{\gamma}
\newcommand{\lmd}{\lambda}
\newcommand{\td}{\tilde}
\newcommand{\vf}{\varphi}
\newcommand{\ad}{\mbox{ad}}
\newcommand{\stl}{\stackrel}
\newcommand{\ol}{\overline}
\newcommand{\es}{\epsilon}
\newcommand{\vsi}{\varsigma}
\newcommand{\ves}{\varepsilon}
\newcommand{\la}{\langle}
\newcommand{\ra}{\rangle}
\newcommand{\vt}{\vartheta}
\newcommand{\wt}{\mbox{wt}\:}
\newcommand{\sym}{\mbox{sym}}
\newcommand{\for}{\mbox{for}}
\newcommand{\rta}{\rightarrow}
\newcommand{\der}{\mbox{Der}}
\newcommand{\mbb}{\mathbb}
\newcommand{\dvg}{\mbox{div}}

\begin{center}{\Large \bf Nongraded Infinite-Dimensional Simple Lie Algebras}\end{center}

\begin{center}{\large Xiaoping Xu}\end{center}

\begin{center}{Department of Mathematics, The Hong Kong University of Science \& Technology}\end{center}
\begin{center}{Clear Water Bay, Kowloon, Hong Kong}\end{center}
\vspace{0.3cm}

\begin{abstract}{Nongraded infinite-dimensional Lie algebras appeared naturally in the theory of Hamiltonian operators, the theory of vertex algebras and their multi-variable analogues. They play important roles in mathematical physics. This survey article is written  based on the author's seminar talks on nongraded infinite-dimensional simple Lie algebras. The key constructional ingredients of our Lie algebras are locally-finite derivations. The structure spaces of some families of these simple Lie algebras can be viewed as analogues of vector bundles 
with Lie algebras as fibers. We also believe that some of our simple Lie algebras could be related to noncommutative geometry.}\end{abstract}

\section{Introduction}

A Lie (super)algebra ${\cal G}$ is called {\it finitely-graded} if ${\cal G}=\bigoplus_{\al\in\G}{\cal G}_{\al}$ is a $\G$-graded vector space for some abelian group $\G$ such that
$$\dim\:{\cal G}_{\al}<\infty,\;\;[{\cal A}_{\al},{\cal A}_{\be}]\subset{\cal A}_{\al+\be}\qquad\for\;\;\al,\be\in \G.\eqno(1.1)$$
In short, we call  finitely-graded as {\it graded}. Graded infinite-dimensional simple Lie (super)algebras have been extensively studied (e.g., cf. [Bl], [CK], [DZ1-5], [K1-4], [K6], [Kn], [L1-2], [Mo1-2], [Mr], [O], [OZ1-4], [S1-3], [Z1], [W]). Our list of references might be partial. An important family of graded infinite-dimensional simple Lie
algebras associated with generalizations of Cartan matrices were introduced by Kac [K1] and Moody [Mr]. Most of the other graded infinite-dimensional simple Lie
algebras are related to derivations and commutative associative algebras. $\mbb{Z}$-graded simple Lie algebras were classified by Kac [K1], Weisfeiler [W] and Mathieu [Mo1]. $\mbb{Z}$-graded-simple Lie algebras were classified by Mathieu [Mo2]. Cheng and Kac [CK] classified $\mbb{Z}$-graded simple Lie superalgebras. 

Vertex  (super)algebras were originally introduced in order to study the moonshine representation of the Monster group (cf. [Bo], [FLM]). Vertex  (super)algebras have become the fundamental algebraic structures in (super)conformal field theory. These algebraic structures have been reformulated by Beilinson and Drinfeld [BD] in terms of algebraic geometry. The Lie (super)algebras associated with  vertex (super)algebras are in general not graded (cf. [Bo], [FLM], [K5], [X3]). They are Lie algebras with one-variable structure.  Kac [K5] introduced the notation of conformal (super)algebra, which is the local structure of a Lie (super)algebra with one-variable structure. Vertex  (super)algebras are generated by conformal (super)algebras. We proved in [X9] that conformal (super)algebras are equivalent to linear Hamiltonian (super)operators. Bakalov, D'Andrea and Kac [BDK] introduced a multi-variable analogue of conformal algebra. 

Our earlier motivation of studying nongraded simple Lie superalgebras was from the study of the Lie superalgebras generated by the Hamiltonian superoperators associated with Novikov-Poisson algebras (cf. [X2], [X5]). The notion of Novikov-Poisson algebra was introduced in [X1] in order to establish a tensor theory for Novikov algebras. More new nongraded Lie (super)algebras appeared in our work [X6] of quadratic conformal superalgebras, which are related to certain Hamiltonian pairs. In general, we are interested in knowing which Lie (super)algebras are simple algebras with one-variable structure. 

A linear transformation $T$ on a vector space $V$ is called {\it locally-finite} if
$$\dim(\mbox{Span}\:\{T^n(u)\mid n\in\mbb{N}\})<\infty\qquad\for\;\;u\in V,\eqno(1.2)$$
where $\mbb{N}$ is the additive semigroup of nonnegative integers. The simple Lie (super)algebras discussed in this article are related to locally-finite derivations of commutative associative algebras. Under certain conditions, they are Lie (super)algebras with multi-variable structures.

Throughout this article, all the vector spaces (algebras) are assumed over a field $\mbb{F}$ with characteristic 0. Most important examples of such $\mbb{F}$ are $\mbb{R}$ the field of real numbers and $\mbb{C}$ the field of complex numbers. We also denote by $\mbb{Z}$ the ring of integers.

\section{Lie Algebras of Witt Type and Divergence-Free}

A {\it derivation} $d$ of a commutative associative algebra ${\cal A}$ is a linear transformation on ${\cal A}$ such that
$$d(uv)=d(u)v+ud(v)\qquad\for\;\;u,v\in{\cal A}.\eqno(2.1)$$
The space $\der\:{\cal A}$ of all the derivations of ${\cal A}$ forms a Lie
algebra with the Lie bracket
$$[d_1,d_2]=d_1d_2-d_2d_1\qquad\for\;\;d_1,d_2\in\der\:{\cal A},\eqno(2.2)$$
where the multiplications on the right-hand side are compositions of linear transformations. For $u\in{\cal A}$ and $d\in \der\:{\cal A}$, we define the operator $ud$ by
$$(ud)(v)=ud(v)\qquad\for\;\;v\in{\cal A}.\eqno(2.3)$$
Then $ud$ is also a derivation of ${\cal A}$. Hence $\der\:{\cal A}$ is a left ${\cal A}$-module.

Given $m,n\in\mbb{Z}$ with $m<n$, we shall
use the
following notation of indices
$$\ol{m,n}=\{m,m+1,m+2,...,n\}\eqno(2.4)$$
throughout this article. We also treat $\ol{m,n}=\emptyset$ when $m>n$.
The {\it classical Witt algebras} is
$${\cal W}(k,\mbb{F})=\der\:\mbb{F}[t_1,t_2,...,t_k]=\sum_{i=1}^k\mbb{F}[t_1,t_2,...,t_k]\ptl_{t_i}.\eqno(2.5)$$
The $\mbb{Z}$-grading ${\cal W}(k,\mbb{F})=\bigoplus_{m\in\mbb{Z}}({\cal W}(k,\mbb{F}))_m$ is defined by
$$({\cal W}(k,\mbb{F}))_m=0\qquad\mbox{if}\;\; -1>m\in\mbb{Z},\eqno(2.6)$$
$$({\cal W}(k,\mbb{F}))_m=\mbox{Span}\{t_1^{n_1}t_2^{n_2}\cdots t_k^{n_k}\ptl_{t_i}\mid i\in\ol{1,k},\;\sum_{j=1}^kn_j=m+1\}\qquad\mbox{if}\;\; -1\leq m\in\mbb{Z}.\eqno(2.7)$$
Then the Witt algebra ${\cal W}(k,\mbb{F})$ is a simple finitely-$\mbb{Z}$-graded Lie algebra.

For a commutative associative algebra ${\cal A}$ and a subspace ${\cal D}$ of derivations, ${\cal A}$ is called ${\cal D}${\it -simple} ({\it derivation-simple} with respect to ${\cal D}$) if there does not exist a subspace ${\cal I}$ of ${\cal A}$ such that ${\cal I}\neq\{0\},{\cal A}$,  and
$$u{\cal I},\;d({\cal I})\subset {\cal I}\qquad\for\;\;u\in{\cal A},\;d\in{\cal D}.\eqno(2.8)$$
When ${\cal D}$ is a Lie subalgebra of $\der\:{\cal A}$, the space
$${\cal W}({\cal A})={\cal A}{\cal D}\eqno(2.9)$$
also form a Lie subalgebra of $\der\:{\cal A}$.
\psp

{\bf Theorem 2.1 (Jordan, [J1-3])}. {\it The Lie algebra ${\cal W}({\cal A})$ is simple if ${\cal A}$ is ${\cal D}$-simple.}
\psp

A deficiency in Jordan's papers is that there are no essentially new examples of simple Lie algebras of Witt type. In [J1-2], Jordan dealt with the case ${\cal D}$ is an abelian subalgebra. The results was repeated in [P1].

For any positive integer $n$, an additive subgroup $G$ of $\mbb{F}^{\:n}$ is called {\it nondegenerate} if $G$ contains an $\mbb{F}$-basis of $\mbb{F}^{\:n}$. Let $\ell_1,\;\ell_2$ and $\ell_3$ be three nonnegative integers such that
$$\ell=\ell_1+\ell_2+\ell_3>0.\eqno(2.10)$$
Take any nondegenerate additive subgroup $\G$ of $\mbb{F}^{\:\ell_2+\ell_3}$ and $\G=\{0\}$ when $\ell_2+\ell_3=0$.  Let ${\cal A}(\ell_1,\ell_2,\ell_3;\G)$ be a free $\mbb{F}[t_1,t_2,...,t_{\ell_1+\ell_2}]$-module with the basis
$$\{x^{\al}\mid \al\in\G\}.\eqno(2.11)$$
Viewing ${\cal A}(\ell_1,\ell_2,\ell_3;\G)$ as a vector space over $\mbb{F}$, we define a commutative associative algebraic operation ``$\cdot$'' on ${\cal A}(\ell_1,\ell_2,\ell_3;\G)$ by
$$(\zeta x^{\al})\cdot (\eta x^{\be})=\zeta\eta x^{\al+\be}\qquad\for\;\;\zeta,\eta\in \mbb{F}[t_1,t_2,...,t_{\ell_1+\ell_2}],\;\al,\be\in\G.\eqno(2.12)$$
Note that  $x^{0}$ is the identity element, which is denoted as $1$ for convenience. When the context is clear, we shall omit the notation ``$\cdot$'' in any associative algebra product.

We define the linear transformations 
$$\{\ptl_{t_1},...,\ptl_{t_{\ell_1+\ell_2}},\ptl^{\ast}_1,...,\ptl^{\ast}_{\ell_2+\ell_3}\}\eqno(2.13)$$
on ${\cal A}(\ell_1,\ell_2,\ell_3;\G)$ by 
$$\ptl_{t_i}(\zeta x^{\al})=\ptl_{t_i}(\zeta)x^{\al},\;\;\ptl^{\ast}_j(\zeta x^{\al})=\al_j \zeta x^{\al}\eqno(2.14)$$
for $\zeta\in \mbb{F}[t_1,t_2,...,t_{\ell_1+\ell_2}]$ and $\al=(\al_1,...,\al_{\ell_2+\ell_3})\in\G.$
 Then $\{\ptl_{t_1},...,\ptl_{t_{\ell_1+\ell_2}},\ptl^{\ast}_1,...,\ptl^{\ast}_{\ell_2+\ell_3}\}$ are mutually commutative derivations of ${\cal A}(\ell_1,\ell_2,\ell_3;\G)$. Set
$$\ptl_i=\ptl_{t_i},\;\;\ptl_{\ell_1+j}=\ptl^{\ast}_j+\ptl_{t_{\ell_1+j}},\;\;\ptl_{\ell_1+\ell_2+l}=\ptl_{\ell_2+l}^{\ast}\eqno(2.15)$$
for $i\in\ol{1,\ell_1},\;j\in\ol{1,\ell_2}$ and $l\in\ol{1,\ell_3}$. Observe that $\{\ptl_i\mid i\in\ol{1,\ell}\}$ is an $\mbb{F}$-linearly independent  mutually commutative set of derivations. 
Let
$${\cal D}=\sum_{i=1}^{\ell}\mbb{F}\ptl_i.\eqno(2.16)$$

In [SXZ], we classified the pairs of a commutative associative algebra ${\cal A}$ with an identity element and a locally-finte finite-dimensional abelian Lie subalgebra ${\cal D}$ of derivations  such that ${\cal A}$ is ${\cal D}$-simple when $\mbb{F}$ is algebraically closed. For neatness, we state here a more restricted result from [SXZ].
\psp

{\bf Theorem 2.2 (Su, Xu and Zhang, [SXZ])}. {\it The algebra ${\cal A}(\ell_1,\ell_2,\ell_3;\G)$ defined in (2.12) is ${\cal D}$-simple with respect to the ${\cal D}$ defined in (2.16). Conversely, suppose that $\mbb{F}$ is algebraically closed,  ${\cal A}$ is a commutative associative algebra  with an identity element and  ${\cal D}$ is a locally-finte finite-dimensional abelian Lie subalgebra ${\cal D}$ of derivations such that ${\cal A}$ is ${\cal D}$-simple and 
$$\bigcap_{\ptl\in {\cal D}}\kn_{_{\ptl}}=\mbb{F},\eqno(2.17)$$
then ${\cal A}$ is isomorphic to ${\cal A}(\ell_1,\ell_2,\ell_3;\G)$ for some $\ell_1,\ell_2,\ell_3\in\mbb{N}$ and addtive subgroup $\G$ of $\mbb{F}^{\:\ell_2+\ell_3}$, and ${\cal D}$ is of the form (2.16).}
\psp

In particular, the Lie algebra 
$${\cal W}(\ell_1,\ell_2,\ell_3;\G)= {\cal A}(\ell_1,\ell_2,\ell_3;\G){\cal D},\eqno(2.18)$$
first given in [X7], is simple, which is not finitely-graded when $\ell_2>0$. The simplest example of such algebra is
$${\cal W}=\sum_{i=1}^k\mbb{F}[t_1^{\pm 1},,...,t_k^{\pm 1}, t_{k+1},...,t_{2k}](t_i\ptl_{t_i}+\ptl_{t_{k+i}}).\eqno(2.19)$$

The following special graded cases were studied by the other people. The algebra ${\cal W}(0,0,\ell_3;\G)$ was introduced by Kawamoto [Kn]. Osborn studied the algebra ${\cal W}(\ell_1,0,\ell_3;\G)$ with
$$\G=\Dlt_1\times\Dlt_1\times\cdots\times\Dlt_{\ell_3},\eqno(2.20)$$
where $\Dlt_i$ are additive subgroups of $\mbb{F}$. The algebra ${\cal W}(\ell_1,0,\ell_3;\G)$ without the condition (2.20) was obtained by Djokovic and Zhao [DZ3]. One can extend the algebra ${\cal W}(\ell_1,\ell_2,\ell_3;\G)$ with ${\cal A}(\ell_1,\ell_2,\ell_3;\G)$ replaced by any of its localization or by adding ${\cal A}(\ell_1,\ell_2,\ell_3;\G)$ with certain powers $u_1^{\al_1}u_2^{\al_2}\cdots u_n^{\al_n}$ for some $0\neq u_i\in {\cal A}(\ell_1,\ell_2,\ell_3;\G)$ and $(\al_1,\al_2,...,\al_n)\in\mbb{F}^{\:n}$. The extended Witt algebras are still simple. Passman [P2] proved a theorem for Lie superalgebras analogous to Jordan's Theorem (Theorem 2.1), where ${\cal D}$ is super-commutative.

Denote by $M_{m\times n}$ the algebra of $m\times n$ matrices with entries in $\mbb{F}$ and by $GL_m$ the group of invertible $m\times m$ matrices. Set
$$G_{\ell_2,\ell_3}=\left\{\left(\begin{array}{cc}A&0_{\ell_2\times\ell_3}\\B& C\end{array}\right)\mid A\in GL_{\ell_2},\;B\in M_{\ell_2\times\ell_3},\;C\in  GL_{\ell_3}\right\},\eqno(2.21)$$
where $0_{\ell_2\times\ell_3}$ is the $\ell_2\times\ell_3$ matrix whose entries are zero. Then $G_{\ell_2,\ell_3}$ forms a subgroup of $GL_{\ell_2+\ell_3}$. Define an action of $G_{\ell_2,\ell_3}$ on $\mbb{F}^{\:\ell_2+\ell_3}$ by
$$ g(\al)=\al g^{-1}\;\;(\mbox{matrix multiplication})\qquad\for\;\;\al\in\mbb{F}^{\:\ell_2+\ell_3},\;g\in G_{\ell_2,\ell_3}.\eqno(2.22)$$
For any nondegenerate additive subgroup $\Upsilon$ of $\mbb{F}^{\:\ell_2+\ell_3}$ and $g\in G_{\ell_2,\ell_2}$, the set
$$g(\Upsilon)=\{g(\al)\mid \al\in \Upsilon\}\eqno(2.23)$$
also forms a nondegenerate additive subgroup of $\mbb{F}^{\:\ell_2+\ell_3}$. Let
$$\Omega_{\ell_2+\ell_3}=\mbox{the set of nondegenerate additive subgroups of}\;\mbb{F}^{\:\ell_2+\ell_3}.\eqno(2.24)$$
We have an action of $G_{\ell_2,\ell_3}$ on $\Omega_{\ell_2+\ell_3}$ by (2.23). Define the moduli space
$${\cal M}_{\ell_2,\ell_3}=\Omega_{\ell_2+\ell_3}/G_{\ell_2,\ell_3},\eqno(2.25)$$
which is the set of $G_{\ell_2,\ell_3}$-orbits in $\Omega_{\ell_2+\ell_3}$. 
\psp

{\bf Theorem 2.3 (Su, Xu and Zhang, [SXZ])}. {\it The Lie algebra} ${\cal W}(\ell_1,\ell_2,\ell_3;\G)$ {\it is isomorphic to the Lie algebra} ${\cal W}(\ell_1',\ell_2',\ell_3';\G')$ {\it if and only if} $(\ell_1,\ell_2,\ell_3)=(\ell_1',\ell_2',\ell_3')$ {\it and there exists an element} $g\in  G_{\ell_2,\ell_3}$ {\it such that}
$g(\G)=\G'$. {\it In particular, there exists a one-to-one correspondence between the set of isomorphism classes of the Lie algebras of the form (2.18) and the following set:}
$$SW=\{(\ell_1,\ell_2,\ell_3,\varpi)\mid (0,0,0)\neq(\ell_1,\ell_2,\ell_3)\in\mbb{N}^{\:3},\;\varpi\in {\cal M}_{\ell_2,\ell_3}\}.\eqno(2.26)$$
{\it In other words, the set} $SW$ {\it is the structure space of the  simple Lie algebras of Witt type in the form (2.18)}. 
\psp

For $\ptl=\sum_{i=1}^nf_i\ptl_{t_i}\in {\cal W}(n,\mbb{F})$ (cf. (2.5)), we define the {\it divergence of} $\ptl$ by
$$\dvg\:\ptl=\sum_{i=1}^n\ptl_{t_i}(f_i).\eqno(2.27)$$
Set 
$${\cal S}(n,\mbb{F})=\{\ptl\in {\cal W}(n,\mbb{F})\mid \dvg\:\ptl =0\}.\eqno(2.28)$$
Then ${\cal S}(n,\mbb{F})$ is the divergence-free Lie subalgebra of ${\cal W}(n,\mbb{F})$. Moreover, ${\cal S}(n,\mbb{F})$ is the subalgebra of ${\cal W}(n,\mbb{F})$ annihilating the volume form; that is,
$${\cal S}(n,\mbb{F})=\{\ptl\in {\cal W}(n,\mbb{F})\mid \ptl(dt_1\wedge dt_2\wedge\cdots\wedge dt_n)=0\}.\eqno(2.29)$$
When $\mbb{F}=\mbb{R}$ the field of real numbers, $e^{\ptl}$ is a volume-preserving diffeomorphism for each $\ptl\in {\cal S}(n,\mbb{F})$. The algebra ${\cal S}(n,\mbb{F})$ is called a {\it Lie algebra of Special type}.

Let us back to our algebra ${\cal W}(\ell_1,\ell_2,\ell_3;\G)$ in (2.18), we define the divergence by
$$\dvg\:\ptl=\sum_{i=1}^{\ell}\ptl_i(u_i)\qquad\for\;\;\ptl=\sum_{i=1}^{\ell}u_i\ptl_i\in {\cal W}(\ell_1,\ell_2,\ell_3;\G),\eqno(2.30)$$
and set
$${\cal S}(\ell_1,\ell_2,\ell_3;\G)=\{\ptl \in {\cal W}(\ell_1,\ell_2,\ell_3;\G)\mid \dvg\:\ptl=0\}.\eqno(2.31)$$
Let $\rho\in\G$ be any element. Then the space
$${\cal S}(\ell_1,\ell_2,\ell_3;\rho,\G)=x^{\rho}{\cal S}(\ell_1,\ell_2,\ell_3;\G)\eqno(2.32)$$
forms a Lie subalgebra of the Lie algebra ${\cal W}(\ell_1,\ell_2,\ell_3;\G)$ (cf. (2.11)).
\psp

{\bf Theorem 2.4. (Xu, [X7])}. {\it The Lie algebra ${\cal S}(\ell_1,\ell_2,\ell_3;\rho,\G)$
 is simple if $\ell_1+\ell_2>0$ or $\rho=0$ or $\ell\geq 3$.  When $\ell_1+\ell_2=0,\;\ell_3=2$ and $\rho\neq 0$,  the derived subalgebra $({\cal S}(\ell_1,\ell_2,\ell_3;\rho,\G))^{(1)}$ is simple and has codimension one in ${\cal S}(\ell_1,\ell_2,\ell_3;\rho,\G)$.} 
\psp

The above theorem was originally proved in [X7] under the condition (2.20) with $\ell_3$ replaced by $\ell_2+\ell_3$. Later we observed that
the condition can be removed by a trick used in [DZ4].  The Lie algebra ${\cal S}(\ell_1,\ell_2,\ell_3;\rho,\G)$ is not finitely-graded if $\ell_2>0$. We also called ${\cal S}(\ell_1,\ell_2,\ell_3;\rho,\G)$ a {\it divergence-free Lie algebra} in [SX1]. The special graded case ${\cal S}(\ell_1,0,\ell_3;\rho,\G)$ was studied by Osborn [O], Djokovic and Zhao [DZ4], and Zhao [Z1].

Define the moduli space
$${\cal M}^{\cal S}_{\ell_2,\ell_3}=(\G\times \Omega_{\ell_2+\ell_3})/G_{\ell_2,\ell_3},\eqno(2.33)$$
where the action of $G_{\ell_2,\ell_3}$ on $\G\times\Omega_{\ell_2+\ell_3}$
is defined by $g(\rho,\Upsilon)=(g(\rho),g(\Upsilon))$
for $(\rho,\Upsilon)\in\G\times\Omega_{\ell_2+\ell_3}$
(cf. (2.21)-(2.24)).
\psp

{\bf Theorem 2.5 (Su and Xu, [SX1])}. {\it The Lie algebras} ${\cal S}(\ell_1,\ell_2,\ell_3;\rho,\G)$
{\it and} ${\cal S}(\ell_1',\ell_2',\ell_3';\rho',\G')$ {\it with} $\ell\geq 3$ {\it are isomorphic if and only if} $(\ell_1,\ell_2,\ell_3)=(\ell_1',\ell_2',\ell_3')$
{\it and there exists an element} $g\in  G_{\ell_2,\ell_3}$ {\it such that}
$g(\G)=\G'$, {\it and} $g(\rho)=\rho'$ {\it if} $\ell_1=0$. {\it In particular, there exists a one-to-one correspondence between the set of isomorphism classes of the Lie algebras of the form (2.32) and the set} $SW$ {\it in (2.26) if} $\ell_1>0$, {\it and between the set of isomorphism classes of the Lie algebras of the form (2.32) and  the following set}:
$$SS=\{(\ell_1,\ell_2,\ell_3,\varpi)\mid (0,0,0)\neq(\ell_1,\ell_2,\ell_3)\in\mbb{N}^{\:3},\;\varpi\in {\cal M}^{\cal S}_{\ell_2,\ell_3}\}\eqno(2.34)$$
{\it if} $\ell_1=0$.
\psp

Bergen and Passman [BP] proved an abstract simplicity theorem for Lie algebras of special type with the assumption that the derivations in ${\cal D}$ are diagonalizable. By our classification theorem of ${\cal D}$-simple algebras in [SXZ] (cf. Theorem 2.2), the Lie algebras of Special type studied in [BP] are essentially the algebras of the form ${\cal S}(0,0,\ell_3;0,\G)$.

\section{Lie Algebras with a Feature of the Block Algebras}

Let $k$ be a positive integer. We shall write an element of $\mbb{F}^{\:k}$ as
$$\al=(\al_1,\al_2,...,\al_k).\eqno(3.1)$$
Denote by ${\cal X}_r$ the $r$th coordinate mapping of $\mbb{F}^{\:k}$ with $r\in\ol{1,k}$, that is,
$${\cal X}_r(\al)=\al_r\qquad\for\;\;\al\in \mbb{F}^{\:k}.\eqno(3.2)$$
For $\lmd\in\mbb{F}$, we let
$$\lmd_{[r]}=(0,...,0,\stl{r}{\lmd},0,...,0)\qquad\for\;\;r\in\ol{1,k}.\eqno(3.3)$$
Take an additive subgroup $\G$ of $\mbb{F}^{\:k}$ and additive semigroups
$${\cal J}_p\in\{\{0\},\mbb{N}\}\qquad \for\;\;p\in\ol{1,k}\eqno(3.4)$$
such that
$${\cal X}_p(\G)+{\cal J}_p\neq\{0\}\qquad \for\;\;p\in\ol{1,k}.\eqno(3.5)$$

Set
$${\cal J}={\cal J}_1\times {\cal J}_2\times\cdots\times {\cal J}_k.\eqno(3.6)$$
Then ${\cal J}$ forms an additive sub-semigroup of $\mbb{F}^{\:k}$. We shall write an element of ${\cal J}$ as
$$\vec i=(i_1,i_2,...,i_k).\eqno(3.7)$$
Let
${\cal A}(\G,{\cal J})$ be the semigroup algebra of $\G\times {\cal J}$ with the canonical basis
$$\{x^{\al,\vec i}\mid (\al,\vec i)\in \G\times {\cal J}\},\eqno(3.8)$$
that is,
$$x^{\al,\vec i}\cdot x^{\be,\vec j}=x^{\al+\be,\vec i+\vec j}\qquad\for\;\;(\al,\vec i),(\be,\vec j)\in \G\times {\cal J}.\eqno(3.9)$$
 Moreover, we have the following mutually commutative derivations $\{\ptl_p\mid p\in\ol{1,k}\}$ defined by
$$\ptl_p(x^{\al,\vec i})=\al_px^{\al,\vec i}+i_px^{\al,\vec i-1_{[p]}}\qquad\for\;\;(\al,\vec i)\in \G\times {\cal J},\eqno(3.10)$$
where we adopt the convention that if a notation is not defined but technically appears in an expression, we always treat it as zero; for instance, $x^{\al,-1_{[1]}}=0$ for any $\al\in\G$.
\psp

{\it Class I}.
\psp

Let $k=2$ and denote ${\cal A}_2={\cal A}(\G,{\cal J})$. We define the following Lie bracket $[\cdot,\cdot]_I$ on ${\cal A}_2$ by
$$[u,v]_I=\ptl_1(u)\ptl_2(v)-\ptl_2(u)\ptl_1(v)+u\ptl_1(v)-\ptl_1(u)v\qquad\for\;\;u,v\in{\cal A}_2.\eqno(3.11)$$
Denote
$$\si_1=(0,1),\qquad\si_2=(0,2).\eqno(3.12)$$
We treat
$$x^{\si_p,\vec 0}=0\qquad\mbox{if}\;\;\si_p\not\in\G.\eqno(3.13)$$
Then $x^{\si_1,\vec 0}$ is a central element of the Lie algebra $({\cal A}_2,[\cdot,\cdot]_I)$. Form the quotient Lie algebra
$${\cal B}_2={\cal A}_2/\mbb{F}x^{\si_1,\vec 0},\eqno(3.14)$$
whose induced Lie bracket is still denoted by $[\cdot,\cdot]_I$.
\psp

{\bf Theorem 3.1}. {\it The Lie algebra $({\cal B}_2,[\cdot,\cdot]_I)$ is simple if ${\cal J}\neq\{\vec 0\}$ or $\si_2\not\in\G$. If ${\cal J}=\{\vec 0\}$ and
$\si_2\in\G$, then the derived subalgebra ${\cal B}_2^{(1)}=[{\cal B},{\cal B}]_I$ is simple and ${\cal B}_2={\cal B}_2^{(1)}\oplus(\mbb{F}x^{\si_2,\vec 0}+\mbb{F}x^{\si_1,\vec 0})$.}
\psp

The above theorem was due to Block [Bl] when ${\cal J}=\{\vec 0\}$ and $\si_2\not\in\G$, due to Djokvic and Zhao [DZ1] when ${\cal J}=\{\vec 0\}$ and $\si_2\in\G$, and due to the author [X4] when ${\cal J}\neq \{\vec 0\}$. The Lie algebra $L^+$ in [Z2] is isomorphic our algebra $({\cal A}_2,[\cdot,\cdot]_I)$ with
$${\cal J}_1=\{0\},\;{\cal J}_2=\mbb{N},\;\;{\cal X}_2(\G)=\{0\}.\eqno(3.15)$$
The general case $L$ in [Z2] can essentially be obtained from the above special case by adding certain powers $(x^{0,1_{[2]}})^{\mu}$ to ${\cal A}_2$ with $\mu\in\mbb{F}$.

If ${\cal X}_p(\G){\cal J}_p\neq \{\vec 0\}$ for $p=1$ or 2, the simple Lie algebra $({\cal B}_2,[\cdot,\cdot]_I)$ is not finitely-graded. The algebra $({\cal A}_2,[\cdot,\cdot]_I)$ appeared in our classification work [X6] of quadratic conformal algebras. The following is a typical example of such algebras.
\psp

{\bf Example 3.1}. Let $m$ be a positive integer. Let ${\cal A}$ be the subalgebra of 
$$\mbb{R}[t_1^{\pm 1/m},t_2^{\pm 1/m},t_3,t_4]\eqno(3.16)$$
 generated by 
$$\{t_1^{\pm 1},t_2^{\pm 1},t_3,t_4,(t_1t_2)^{1/m}\}.\eqno(3.17)$$
In particular, ${\cal A}=\mbb{R}[t_1^{\pm 1},t_2^{\pm 1},t_3,t_4]$ when $m=1$. The space ${\cal B}={\cal A}/\mbb{R}t_2$ and its Lie bracket is induced by 
$$[f,g]_I=(t_1f_{t_1}+f_{t_3})(t_2g_{t_2}-g_{t_4}-g)+(f-t_2f_{t_2}-f_{t_4})(t_1g_{t_1}+g_{t_3})\eqno(3.18)$$
for $f,g\in{\cal A}$.
\psp

{\it Class II}.
\psp

Let $k=4$ and denote ${\cal A}_4={\cal A}(\G,{\cal J})$. Assume 
$$\mbb{Z}^4\subset \G.\eqno(3.19)$$
Take
$$0\neq (\kappa_1,\kappa_2,0,0),\;0\neq (0,0,\kappa_3,\kappa_4)\in\G\eqno(3.20)$$
and set
$$\kappa=(\kappa_1,\kappa_2,\kappa_3,\kappa_4).\eqno(3.21)$$
The following Lie bracket $[\cdot,\cdot]_{II}$ on ${\cal A}_4$ was obtained in [X4] by the author in studying quadratic conformal algebras (cf. [X6]):
$$[u,v]_{II}=x^{\kappa,\vec 0}(\ptl_1(u)\ptl_2(v)-\ptl_2(u)\ptl_1(v))+(\ptl_3+\kappa_3)(u)(\ptl_4+\kappa_4)(u)-(\ptl_4+\kappa_4)(u)(\ptl_3+\kappa_3)(u)\eqno(3.22)$$
for $u,v\in{\cal A}_4$. Set 
$$\si=(0,0,-\kappa_3,-\kappa_4),\;\;\rho=(\kappa_1,\kappa_2,-2\kappa_3,-2\kappa_4).\eqno(3.23)$$
Then $x^{\si,\vec 0}$ is a central element of $({\cal A}_4,[\cdot,\cdot]_{II})$. Form the quotient Lie algebra
$${\cal B}_4={\cal A}_4/\mbb{F}x^{\si_1,\vec 0},\eqno(3.24)$$
whose induced Lie bracket is still denoted by $[\cdot,\cdot]_{II}$.
\psp

{\bf Theorem 3.2 (Xu, [X4])}. {\it The Lie algebra $({\cal B}_4,[\cdot,\cdot]_{II})$ is simple if ${\cal J}\neq\{\vec 0\}$. If ${\cal J}=\{\vec 0\}$, then the derived subalgebra ${\cal B}_4^{(1)}=[{\cal B},{\cal B}]_{II}$ is simple and ${\cal B}_4={\cal B}_4^{(1)}\oplus(\mbb{F}x^{\rho,\vec 0}+\mbb{F}x^{\si,\vec 0})$.}
\psp

The Lie bracket $[\cdot,\cdot]_{II}$ in (3.22) may not look natural from pure Lie algebra point of view. However, it is quite natural from quadratic conformal algebra point of view. The finding of this family of simple Lie algebras shows that the notion of conformal algebra is useful in studying simple Lie algebras. The Lie algebra $({\cal B}_4,[\cdot,\cdot]_{II})$ is not finitely-graded.
\psp

{\bf Example 3.2}. Let ${\cal A}$ be the subalgebra of 
$$\mbb{R}[t_p^{\pm 1/m},t_{4+p}\mid p\in\ol{1,4}]\eqno(3.25)$$
 generated by 
$$\{t_p^{\pm 1},t_{4+p},(t_1t_2t_3t_4)^{1/m}\mid p\in\ol{1,4}\}.\eqno(3.26)$$
When $m=1$,
$${\cal A}=\mbb{R}[t_p^{\pm 1},t_{4+p}\mid p\in \ol{1,4}].\eqno(3.27)$$
Take a nonzero integer $n$. The space ${\cal B}={\cal A}/\mbb{R}(t_3t_4)^{-n}$ and its Lie bracket is induced by
\begin{eqnarray*}[f,g]_{II}&=&(t_1t_2t_3t_4)^n[(t_1f_{t_1}+f_{t_5})(t_2g_{t_2}+g_{t_6})-(t_2f_{t_2}+f_{t_6})(t_1g_{t_1}+g_{t_5})]+(t_3f_{t_3}+\\ & &+f_{t_7}+nf)(t_4g_{t_4}+g_{t_8}+ng)-(t_4f_{t_4}+f_{t_8}+nf)(t_3g_{t_3}+g_{t_7}+ng)\hspace{1.9cm}(3.28)\end{eqnarray*}
for $f,g\in{\cal A}$. 
\psp

{\it Class III}. 
\psp

The third class of algebras are supersymmetric extensions of the algebras of Class I.
Set
$$\td{\cal A}={\cal A}_2\times {\cal A}_2=\td{\cal A}_0\oplus \td{\cal A}_1\eqno(3.29)$$
with
$$\td{\cal A}_0=({\cal A}_2,0),\qquad \td{\cal A}_1=(0,{\cal A}_2).\eqno(3.30)$$
Moreover, we denote
$$u_{(0)}=(u,0),\;\;u_{(1)}=(0,u)\eqno(3.31)$$
for $u\in {\cal A}_2.$

Take an element $\kappa=(\kappa_1,\kappa_2)\in\G$.  We have the following Lie superbracket $[\cdot,\cdot]$ on $\td{\cal A}$ defined by
$$[u_{(0)},v_{(0)}]=[\ptl_1(u)\ptl_2(v)-\ptl_1(v)\ptl_2(u)+u\ptl_1(v)-\ptl_1(u)v]_{(0)},\eqno(3.32)$$
$$[u_{(0)},v_{(1)}]=[\ptl_1(u)(\ptl_2(v)+(\kappa_2-1)v/2)+(u-\ptl_2(u))(\ptl_1(v)+\kappa_1v/2)]_{(1)},\eqno(3.33)$$
$$[u_{(1)},v_{(1)}]=(x^{\kappa,\vec{0}}uv)_{(0)}\eqno(3.34)$$
for $u,v\in{\cal A}_2$.  The Lie superalgebras generated by Hamiltonian superoperator of supervaribles (cf. [X5]) associated with certain Novikov-Poisson algebras are of the above form  $(\td{\cal A}, [\cdot,\cdot])$ (cf. [X2]). Set
$$\td{\cal B}=\td{\cal A}_0+\td{\cal B}_1\qquad\mbox{with}\;\;\td{\cal B}_1=[\td{\cal A}_0,{\cal A}_1].\eqno(3.35)$$
Then $\td{\cal B}$ form a Lie sub-superalgebra of $\td{\cal A}$. Moreover,  $(x^{\si_1,\vec 0})_{(0)}$ a central element of $\td{\cal B}$ (cf. (3.12), (3.13)). Form the quotient Lie superalgebra
$$\td{\cal C}=\td{\cal B}/\mbb{F}(x^{\si_1,\vec 0})_{(0)},\eqno(3.36)$$
whose induced Lie superbracket is still denoted by $[\cdot,\cdot]$. Now we allow ${\cal X}_2(\G)={\cal J}_2=\{0\}$ (cf. (3.2), (3.4), (3.5)).
\psp

{\bf Theorem 3.3 (Xu, [X4])}. {\it The Lie superalgebra} $(\td{\cal C},[\cdot,\cdot])$ is simple. Moreover, $\td{\cal B}_1=\td{\cal A}_1$ when ${\cal J}\neq\{\vec 0\}$, and
$$\td{\cal B}_1=\mbox{Span}\:\{(x^{\al,\vec 0})_{(1)}\mid ((0,3)-\kappa)/2\neq \al\in\G\}\eqno(3.37)$$
if ${\cal J}=\{\vec 0\}$. 
\psp

When
$$ {\cal J}=\{\vec 0\},\;\;\G=(\mbb{Z},0),\;\;\kappa=(1,0),\eqno(3.38)$$
the Lie superalgebra $\td{\cal C}=\td{\cal A}$ is the well-known centerless super Virasoro algebra.
Kantor [Ki] constructed Jordan and Lie superalgebras from Poisson algebras. His construction is now known as {\it Kantor's double processing}.  Our above construction of Lie superalgebras from the algebras of Class I is another double processing. Note that $({\cal A}_2,\cdot,[\cdot,\cdot]_I)$ does not form a Poisson algebra (cf. (3.9), (3.11)). The Lie superalgebra $\td{\cal C}$ is not finitely generated when ${\cal X}_p(\G){\cal J}_p\neq \{\vec 0\}$ for $p=1$ or 2.
\psp

{\bf Example 3.3}. Let $k$ be a positive integer and let ${\cal A}$ be the subalgebra of 
$$\mbb{R}[t_1^{\pm 1/k},t_2^{\pm 1/k},t_3,t_4]\eqno(3.39)$$
 generated by 
$$\{t_1^{\pm 1},t_2^{\pm 1},t_3,t_4,(t_1t_2)^{1/k}\}.\eqno(3.40)$$
In particular, ${\cal A}=\mbb{R}[t_1^{\pm 1},t_2^{\pm 1},t_3,t_4]$ when $k=1$. Suppose that $m$ and $n$ are two fixed integers. The space 
$$\bar{\cal B}=({\cal A}\times {\cal A})/\mbb{R}(t_2)_{(0)}\eqno(3.41)$$
with its Lie superbracket induced by
$$[f_{(0)},g_{(0)}]=[(t_1f_{t_1}+f_{t_3})(t_2g_{t_2}+g_{t_4}-g)+(f-t_2f_{t_2}-f_{t_4})(t_1g_{t_1}+g_{t_3})]_{(0)},\eqno(3.42)$$
$$[f_{(0)},g_{(1)}]=[(t_1f_{t_1}+f_{t_3})(t_2g_{t_2}+g_{t_4}+(n-1)g/2)+(f-t_2f_{t_2}-f_{t_4})(t_1g_{t_1}+g_{t_3}+mg/2)]_{(1)},\eqno(3.43)$$
$$[f_{(1)},g_{(1)}]=((t_1^mt_2^n)fg)_{(0)},\eqno(3.44)$$
for $f,g\in{\cal A}$ (cf. (3.31)).

\section{Lie Algebras of Hamiltonian Type and Contact Type}

Let $k$ be a nonnegative integer. We shall write an element of $\mbb{F}^{\:2k}$ as
$$\al=(\al_1,\al_2,...,\al_{2k}).\eqno(4.1)$$
Denote by ${\cal X}_r$ the $r$th coordinate mapping of $\mbb{F}^{\:2k}$ with $r\in\ol{1,2k}$, that is,
$${\cal X}_r(\al)=\al_r\qquad\for\;\;\al\in \mbb{F}^{\:2k}.\eqno(4.2)$$
For $\lmd\in\mbb{F}$, we let
$$\lmd_{[r]}=(0,...,0,\stl{r}{\lmd},0,...,0)\qquad\for\;\;r\in\ol{1,2k}.\eqno(4.3)$$
Take an additive subgroup $\G_1$ of $\mbb{F}^{\:2k}$ and additive semigroups
$${\cal J}_p\in\{\{0\},\mbb{N}\}\qquad \for\;\;p\in\ol{1,2k}\eqno(4.4)$$
such that
$${\cal X}_p(\G)+{\cal J}_p\neq\{0\}\qquad \for\;\;p\in\ol{1,2k}.\eqno(4.5)$$

Set
$${\cal J}={\cal J}_1\times {\cal J}_2\times\cdots\times {\cal J}_{2k}.\eqno(4.6)$$
Then ${\cal J}$ forms an additive sub-semigroup of $\mbb{F}^{\:2k}$. We shall write an element of ${\cal J}$ as
$$\vec i=(i_1,i_2,...,i_{2k}).\eqno(4.7)$$
Let $\G_0$ be a torsion-free abelian group. Define
$$\G=\G_0\oplus\G_1\;\;\;\mbox{a direct sum of abelian groups}.\eqno(4.8)$$
We write an element of $\G$ as
$$\vec \al=\al_0+\al\qquad\mbox{with}\;\;\al_0\in\G_0,\;\al\in\G_1.\eqno(4.9)$$

Let ${\cal A}$ be the semigroup algebra of $\G\times {\cal J}$ with the canonical basis
$$\{x^{\vec\al,\vec i}\mid (\vec\al,\vec i)\in \G\times {\cal J}\},\eqno(4.10)$$
that is,
$$x^{\vec\al,\vec i}\cdot x^{\vec\be,\vec j}=x^{\vec\al+\vec \be,\vec i+\vec j}\qquad\for\;\;(\vec\al,\vec i),(\vec\be,\vec j)\in \G\times {\cal J}.\eqno(4.11)$$
The identity element of ${\cal A}(\G,{\cal J})$ is $x^{\vec 0,\vec 0}$, which will be simply denoted as $1_{\cal A}$. Moreover, we have the following mutually commutative derivations $\{\ptl_p\mid p\in\ol{1,2k}\}$ defined by
$$\ptl_p(x^{\vec\al,\vec i})=\al_px^{\vec \al,\vec i}+i_px^{\vec\al,\vec i-1_{[p]}}\qquad\for\;\;(\vec\al,\vec i)\in \G\times {\cal J},\eqno(4.12)$$
where we adopt the convention that if a notation is not defined but technically appears in an expression, we always treat it as zero; for instance, $x^{\vec\al,-1_{[1]}}=0$ for any $\vec\al\in\G$.

Take $k_1\in\ol{0,k}$ such that
$${\cal X}_p(\G_1)\neq\{0\}\;\;\mbox{or}\;\;{\cal X}_{k+p}(\G_1)\neq\{0\}\qquad\for\;\;p\in\ol{1,k_1}.\eqno(4.13)$$
Set
$$\mho=\{p,k+p\mid p\in\ol{1,k},\;{\cal X}_p(\G_1){\cal X}_{k+p}(\G_1)\neq\{0\}\}.\eqno(4.14)$$
We also view $\G$ as a $\mbb{Z}$-module. Take a skew-symmetric $\mbb{Z}$-bilinear $\phi(\cdot,\cdot):\G\times\G\rta\mbb{F}$ such that
$$\{\al_0\in\G_0\mid\phi(\al_0,\vec\be)=0\;\for\;\vec\be=\be_0+\be\in\G\;\mbox{with}\;{\cal X}_q(\be)=0\;\mbox{if}\;q\in\mho\}=\{0\}.\eqno(4.15)$$
Moreover, we assume
$$1_{[p]}\in\G_1\;\;\mbox{if}\;\;p\in\ol{k_1+1,k}\bigcup\ol{k+k_1+1,2k}\;\;\mbox{and}\;\;{\cal X}_p(\G_1)\neq\{0\},\eqno(4.16)$$
$$1_{[p]}\in\rd_{\phi}\;\;\mbox{if}\;\;p\in\ol{1,k_1}\bigcup\ol{k+1,k+k_1}\;\;\mbox{and}\;\;{\cal X}_p(\G_1)\neq\{0\}.\eqno(4.17)$$

Choose
$$0\neq \si_p\in\rd_{\phi}\bigcap(\mbb{F}1_{[p]}+\mbb{F}1_{[k+p]})\;\;\for\;\;p\in\ol{1,k_1},\;\;\;\si_q=0\;\;\for\;q\in\ol{k_1+1,k},\eqno(4.18)$$
and set
$$\si=\sum_{p=1}^k\si_p.\eqno(4.19)$$
For any $\vec \al\in\G$, we set
$${\cal A}_{\vec\al}=\sum_{\vec i\in{\cal J}}\mbb{F}x^{\vec\al,\vec i}.\eqno(4.20)$$
Then we have the following Lie bracket $[\cdot,\cdot]_H$ on ${\cal A}$:
\begin{eqnarray*}\hspace{3cm}[u,v]_H&=&\sum_{p=1}^kx^{\si_p,\vec 0}(\ptl_p(u)\ptl_{k+p}(v)-\ptl_{k+p}(u)\ptl_p(v))\\ & & +(\phi(\vec\al,\vec\be)-\sum_{p=k_1+1}^k(\al_p\be_{k+p}-\al_{k+p}\be_p))uv\hspace{3cm}(4.21)\end{eqnarray*}
Obviously, $1_{\cal A}$ is a central element of $({\cal A},[\cdot.\cdot]_H)$. Form the quotient Lie algebra
$${\cal H}={\cal A}/\mbb{F}1_{\cal A},\eqno(4.22)$$
whose induced Lie bracket is still denoted by $[\cdot,\cdot]_H$. The algebra ${\cal H}$ is called a {\it Lie algebra of Hamiltonian type}.
\psp

{\bf Theorem 4.1 (Xu, [X7])}. {\it The Lie algebra $({\cal H},[\cdot,\cdot]_H)$ is simple if ${\cal J}\neq\{\vec 0\}$ or $k_1=0$; If ${\cal J}=\{\vec 0\}$ and $k_1>0$, then the derived subalgebra ${\cal H}^{(1)}=[{\cal H},{\cal H}]_H$ is simple and ${\cal H}={\cal H}^{(1)}\oplus(\mbb{F}x^{\si,\vec 0}+\mbb{F}1_{\cal A})$.}
\psp

A classical Lie algebra of Hamiltonian type is isomorphic to the above ${\cal H}$ when
$$\G_0=\{0\},\;\;k_1=0,\;\;\G_1=\{0\}.\eqno(4.23)$$
In this case, ${\cal A}\cong \mbb{F}[t_1,t_2,...,t_{2k}]$. Osborn [O] studied the case
$$\G_0=\{0\},\;\;{\cal X}_p(\G_1)1_{[p]}\subset \G_1,\;\; {\cal X}_p(\G_1)=\{0\}\;\;\mbox{or}\;\;{\cal J}_p=\{0\}\eqno(4.24)$$
for $p\in\ol{1,2k}$. Osborn and Zhao [OZ2] dealt with the case
$${\cal J}=\{\vec 0\},\;\;\phi(\al_0+\al,\be+\be_0)=\phi'(\al_0,\be_0)+\sum_{p=k_1+1}^k(\al_p\be_{k+p}-\al_{k+p}\be_p)\eqno(4.25)$$
for $\al_0,\be_0\in\G_0$ and $\al,\be\in\G_1,$ where $\phi'(\cdot,\cdot):\G_0\times\G_0\rta\mbb{F}$ is a skew-symmetric $\mbb{Z}$-bilinear form. In fact, ${\cal H}$ is not finitely-graded if ${\cal X}_p(\G){\cal J}_p\neq \{\vec 0\}$ for some $p\in\ol{1,2k}$.

It is straightforward to verify
$$[u,vw]_H=[u,v]_Hw+v[u,w]_H\qquad\for\;\;u,v,w\in{\cal A}.\eqno(4.26)$$
Thus $({\cal A},\cdot,[\cdot,\cdot]_H)$ forms a Poisson algebra (cf. (4.11) and (4.21)). In [SX2], we completely determined the isomorphism classes of Poisson algebras of the above form $({\cal A},\cdot,[\cdot,\cdot]_H)$.
\psp

Next we shall construct simple Lie algebras of contact type. Let $k$ be a positive integer. Take an additive subgroup $\G_0$ of $\mbb{F}$, an additive subgroup $\G_1$ of $\mbb{F}^{\:2k}$, and additive semigroups
$${\cal J}_p\in\{\{0\},\mbb{N}\}\qquad \for\;\;p\in\ol{0,2k}\eqno(4.27)$$
such that
$$\G_0+{\cal J}_0\neq\{0\},\;\;{\cal X}_p(\G)+{\cal J}_p\neq\{0\}\qquad \for\;\;p\in\ol{1,2k}.\eqno(4.28)$$

Set
$$\G=\G_0\times\G_1\eqno(4.29)$$
and
$${\cal J}={\cal J}_0\times {\cal J}_1\times\cdots\times {\cal J}_{2k}.\eqno(4.39)$$
Then $\G$ forms a subgroup of $\mbb{F}^{\:2k+1}$ and ${\cal J}$ forms an additive sub-semigroup of $\mbb{F}^{\:2k+1}$. We shall write an element of $\G$ as
$$\al=(\al_0,\al_1,...,\al_{2k})\eqno(4.31)$$
and an element of ${\cal J}$ as
$$\vec i=(i_0,i_1,...,i_{2k}).\eqno(4.32)$$

Let ${\cal A}$ be the semigroup algebra of $\G\times {\cal J}$ with the canonical basis
$$\{x^{\al,\vec i}\mid (\al,\vec i)\in \G\times {\cal J}\},\eqno(4.33)$$
that is,
$$x^{\al,\vec i}\cdot x^{\be,\vec j}=x^{\al+\vec \be,\vec i+\vec j}\qquad\for\;\;(\al,\vec i),(\vec\be,\vec j)\in \G\times {\cal J}.\eqno(4.34)$$
 For $\lmd\in\mbb{F}$, we let
$$\lmd_{[r]}=(0,...,0,\stl{r}{\lmd},0,...,0)\qquad\for\;\;r\in\ol{0,2k}.\eqno(4.35)$$
Moreover, we have the following mutually commutative derivations $\{\ptl_p\mid p\in\ol{0,2k}\}$ defined by
$$\ptl_p(x^{\al,\vec i})=\al_px^{\al,\vec i}+i_px^{\al,\vec i-1_{[p]}}\qquad\for\;\;(\al,\vec i)\in \G\times {\cal J}.\eqno(4.36)$$

We assume that
$$\si_p=-1_{[p]}\in \G\;\;\mbox{if}\;\;{\cal X}_p(\G_1)\neq\{0\},\;{\cal X}_{k+p}(\G_1)=\{0\},\eqno(4.37)$$
$$\si_p=-1_{[k+p]}\in \G\;\;\mbox{if}\;\;{\cal X}_p(\G_1)=\{0\},\;{\cal X}_{k+p}(\G_1)\neq\{0\},\eqno(4.38)$$
$$\si_p=-1_{[p]}-1_{[k+p]}\in\G,\;\;\mbb{F}1_{[p]}\bigcap\G\neq\{0\},\;\mbb{F}1_{[k+p]}\bigcap\G\neq\{0\}\;\;\mbox{if}\;\;{\cal X}_p(\G_1){\cal X}_{k+p}(\G_1)\neq \{0\},\eqno(4.39)$$
for $p\in\ol{1,k}$. 

Set
$$\mho_1=\{p\in\ol{1,2k}\mid {\cal X}_p(\G_1)\neq\{0\}\},\qquad\mho_2=\ol{1,2k}\setminus\mho_1.\eqno(4.40)$$
Assumption (4.28) implies ${\cal J}_p=\mbb{N}$ for $p\in\mho_2$. We define another derivation $\ptl$ on ${\cal A}$ by 
$$\ptl(x^{\al,\vec i})=(\sum_{p\in\mho_1}\al_p+\sum_{q\in\mho_2}i_q)x^{\al,\vec i}\qquad\for\;\;(\al,\vec i)\in\G\times{\cal J}.\eqno(4.41)$$
Take any $\si_0\in\G_01_{[0]}$ and let
$$\si_p=0\qquad\for\;\;p\in\ol{1,k}\setminus\mho_1.\eqno(4.42)$$
Then we have the following Lie bracket $[\cdot,\cdot]_K$ on ${\cal A}$:
$$[u,v]_K=x^{\si_0,\vec 0}[\ptl_0(u)(2-\ptl)(v)-(2-\ptl)(u)\ptl_0(v)]+\sum_{p=1}^kx^{\si_p,\vec 0}(\ptl_p(u)\ptl_{k+p}(v)-\ptl_{k+p}(u)\ptl_p(v))\eqno(4.43)$$
for $u,v\in {\cal A}$. We call the algebra $({\cal A},[\cdot,\cdot]_K)$ a {\it Lie algebra of Contact type}.
\psp

{\bf Theorem 4.2 (Xu, [X7])}. {\it The Lie algebra $({\cal A},[\cdot,\cdot]_K)$ is simple}.
\psp

The Lie algebra $({\cal A},[\cdot,\cdot]_K)$ is not finitely graded if $\G_0{\cal J}_0\neq\{0\}$ or ${\cal X}_p(\G_1){\cal J}_p\neq \{0\}$ for some $p\in\ol{1,2k}$. A classical Contact Lie algebra is isomorphic to $({\cal A},[\cdot,\cdot]_K)$ when $\G=\{0\}$. Osborn [O] studied the algebra $({\cal A},[\cdot,\cdot]_K)$ when $\G_0{\cal J}_0=\{0\}$ and
$${\cal X}_p(\G_1)1_{[p]}\in\G_1,\;\;{\cal X}_p(\G_1){\cal J}_p=\{0\}\qquad\for\;\; p\in\ol{1,2k}.\eqno(4.44)$$
Osborn and Zhao [OZ1] dealt with the case ${\cal J}=\{\vec 0\}$.

\section{Classical Lie Algebras over Generalized Weyl Algebras}

For an associative algebra ${\cal B}$ (which may not have an identity element), the associated
Lie bracket is defined by
$$[u,v]=uv-vu\qquad\for\;\;u,v\in{\cal B}.\eqno(5.1)$$
The {\it Weyl algebra} $\mbb{A}_k$ of rank $k$ is an associative
algebra
generated by $\{t_i,\ptl_i\mid i\in\ol{1,k}\}$ with the defining relations:
$$[t_i,t_j]=0=[\ptl_i,\ptl_j],\,\,\ptl_it_j-t_j\ptl_i=\dlt_{i,j}1_{ \mbb{A}_k },\qquad\for\;\;i,j\in\ol{1,k}.\eqno(5.2)$$

Let us go back to the settings (2.10)-(2.16). For convenience, we denote
$${\cal A}={\cal A}(\ell_1,\ell_2,\ell_3;\G),\eqno(5.3)$$
the commutative associative algebra defined in (2.11) and (2.12). Recall that the derivation subalgebra ${\cal D}$ defined by (2.13)-(2.16).
Now we define
$$\mbb{A}=\sum_{n=0}^{\infty}{\cal A}{\cal D}^n\subset \mbox{End}\;{\cal A}.\eqno(5.4)$$
Then $\mbb{A}$ forms an associative subalgebra of $\mbox{End}\;{\cal A}$. In fact, $\mbb{A}$ is isomorphic to the Weyl algebra
$\mbb{A}_{\ell}$ when $\ell_2=\ell_3=0$. We call the algebra $\mbb{A}$ a {\it generalized Weyl algebra}. 
\psp

For positive integer $n$, we denote by $M_{n\times n}(\mbb{A})$  the algebra of $n\times n$ matrices with entries in $\mbb{A}$. Fix a positive integer $k$. Let 
$$\{\mbb{B}_l^{(1)},\mbb{B}_l^{(2)},...,\mbb{B}^{(k)}_l\}\eqno(5.5)$$
be a set of nonzero left ideals of $\mbb{A}$ and let
$$\{\mbb{B}_r^{(1)},\mbb{B}_r^{(2)},...,\mbb{B}_r^{(k)}\}.\eqno(5.6)$$

 For notational convenience, we set
$$\vec\mbb{B}_l=(\mbb{B}_l^{(1)},\mbb{B}_l^{(2)},...,\mbb{B}^{(k)}_l),\;\;\vec\mbb{B}_r=(\mbb{B}_r^{(1)},\mbb{B}_r^{(2)},...,\mbb{B}^{(k)}_r).\eqno(5.7)$$
We define
$$gl(\vec\mbb{B}_l,\vec\mbb{B}_r)=\left\{\left(\begin{array}{cccc}a_{1,1}& a_{1,2}&\cdots& a_{1,k}\\ a_{2,1}&a_{2,2}&\cdots& a_{2,k}\\ \vdots&\vdots& &\vdots\\  a_{k,1}&a_{k,2}&\cdots& a_{k,k}\end{array}\right)\mid a_{i,j}\in\mbb{B}_r^{(i)}\mbb{B}_l^{(j)}\right\}.\eqno(5.8)$$
Then $gl(\vec\mbb{B}_l,\vec\mbb{B}_r)$ forms a Lie subalgebra of $M_{k\times k}(\mbb{A})$  with the Lie bracket defined in (5.1). Set
$$\vec\mbb{A}=(\mbb{A},\mbb{A},...,\mbb{A}).\eqno(5.9)$$
When $(\vec\mbb{B}_l,\vec\mbb{B}_r)=\vec \mbb{A}$, $gl(\vec\mbb{B}_l,\vec\mbb{B}_r)=M_{k\times k}(\mbb{A})$ has a one-dimensional center $\mbb{F}I_k$, where $I_k$ is the $k\times k$ identity matrix. We define the quotient Lie algebra
$$sl(k,\mbb{A})=M_{k\times k}(\mbb{A})/\mbb{F}I_k.\eqno(5.10)$$

{\bf Theorem 5.1 (Benkart, Xu and Zhao, [BXZ])}. {\it The Lie algebra $gl(\vec\mbb{B}_l,\vec\mbb{B}_r)$ with $(\vec\mbb{B}_l,\vec\mbb{B}_r)\neq \vec\mbb{A}$  and the Lie algebra $sl(k,\mbb{A})$ are simple.}
\psp

Pick
$$\ell'\in \ol{0,\ell_1}.\eqno(5.11)$$
Set
$${\cal A}'=\sum_{\al\in \G}\mbb{F}[t_{\ell'+1},...,t_{\ell_1+\ell_2}]x^{\al},\;\;\mbb{A}'={\cal A}'(\mbb{F}[\ptl_{\ell'+1},...,\ptl_{\ell}]).\eqno(5.12)$$
Then $\mbb{A}'$ forms a  subalgebra of $\mbb{A}$. Define
$$\bar{t}^{\vec m}=t_1^{m_1}t_2^{m_2}\cdots t_{\ell'}^{m_{\ell'}},\;\;\bar{\ptl}^{\vec m}=\ptl_1^{m_1}\ptl_2^{m_2}\cdots \ptl_{\ell'}^{m_{\ell'}}\qquad\for\;\;\vec m=(m_1,...,m_{\ell'})\in\mbb{N}^{\:\ell'},\eqno(5.13)$$
$$\td{\ptl}^{\vec n}=\ptl_{\ell'+1}^{n_{\ell'+1}}\ptl_{\ell'+2}^{n_{\ell'+2}}\cdots \ptl_{\ell}^{n_{\ell}},\;\;(-\td\ptl)^{\vec n}=(-\ptl_{\ell'+1})^{n_{\ell'+1}}(-\ptl_{\ell'+2})^{n_{\ell'+2}}\cdots (-\ptl_{\ell})^{n_{\ell}}\eqno(5.14)$$
for $\vec n=(n_{\ell'+1},...,n_{\ell})\in\mbb{N}^{\:\ell-\ell'}$. Set
$$\mbb{A}''=\sum_{\vec m',\vec m''\in\mbb{N}^{\:\ell'}}\mbb{F}\bar{t}^{\vec m'}\bar{\ptl}^{\vec m''}.\eqno(5.15)$$
Observe that $\mbb{A}''$ is a subalgebra of $\mbb{A}$. Moreover, $\mbb{A}'$ and $\mbb{A}''$ are mutually commutative subalgebras and
$$\mbb{A}=\mbb{A}''\mbb{A}'.\eqno(5.16)$$
We define $\tau\in\mbox{End}\:\mbb{A}$ by
$$\tau(u\bar{t}^{\vec m'}\bar{\ptl}^{\vec m''}\td{\ptl}^{\vec n})=(-\td{\ptl})^{\vec n}\cdot u\bar{t}^{\vec m''}\bar{\ptl}^{\vec m'}\qquad\for\;\;\vec m',\vec m''\in\mbb{N}^{\:\ell'},\;\vec n\in\mbb{N}^{\:\ell-\ell'},\;u\in{\cal A}'.\eqno(5.17)$$
It was verified in [BXZ] that the map $\tau$ {\it is an involutive anti-automorphism of} $\mbb{A}$, {\it that is},
$$\tau^2=1,\;\;\tau(1)=1,\;\;\tau(uv)=\tau(v)\tau(u)\qquad\mbox{\it for}\;\;u,v\in\mbb{A}.\eqno(5.18)$$

Fix a positive integer $k$. Let $\{\mbb{B}_l^{(1)},\mbb{B}_l^{(2)},...,\mbb{B}^{(k)}_l\}$ be a set of nonzero left ideals of $\mbb{A}$ and let $\{\mbb{B}_r^{(1)},\mbb{B}_r^{(2)},...,\mbb{B}_r^{(k)}\}$ be a set of nonzero right ideals of $\mbb{A}$.
Set
$$\hat\mbb{A}=\bigoplus_{i,j=1}^k\mbb{B}_r^{(i)}\mbb{B}_l^{(j)}\eqno(5.19)$$
as a direct sum of $k^2$ vector spaces. Let $\rho\in \mbox{End}\:\hat\mbb{A}$ such 
$$\rho^2=\mbox{Id}_{\hat{A}},\;\;\rho(\mbb{B}_r^{(i)}\mbb{B}_l^{(j)})=\mbb{B}_r^{(j)}\mbb{B}_l^{(i)}\qquad\for\;\;i,j\in\ol{1,k}\eqno(5.20)$$
and
$$\rho(ab)=\rho(b)\rho(a)\qquad\for\;\;i,j,p\in\ol{1,k},\;a\in \mbb{B}_r^{(i)}\mbb{B}_l^{(j)},\;b\in \mbb{B}_r^{(j)}\mbb{B}_l^{(p)}.\eqno(5.21)$$
The followings are examples of such a linear map $\rho$.
\psp

{\bf Example 5.1}. (1) We first take  a set $\{\mbb{B}_l^{(1)},\mbb{B}_l^{(2)},...,\mbb{B}^{(k)}_l\}$  of nonzero left ideals of $\mbb{A}$, and then take
$$\mbb{B}_r^{(j)}=\tau(\mbb{B}_l^{(j)})\qquad\for\;\;j\in\ol{1,k}.\eqno(5.22)$$
It is easy to see that $\{\mbb{B}_r^{(1)},\mbb{B}_r^{(2)},...,\mbb{B}^{(k)}_r\}$ is a set of nonzero right ideals of $\mbb{A}$. Moreover,
$$\tau(\mbb{B}_r^{(i)}\mbb{B}_l^{(j)})=\mbb{B}_r^{(j)}\mbb{B}_l^{(i)}\qquad\for\;\;i,j\in\ol{1,k}\eqno(5.23)$$
by (5.18). So we can take $\rho\in \mbox{End}\:\hat\mbb{A}$ determined by
$$\rho(u)=\tau(u)\qquad \for\;\;u\in \mbb{B}_r^{(i)}\mbb{B}_l^{(j)},\;i,j\in\ol{1,k}.\eqno(5.24)$$

(2) Take
$$\vec m_1,\vec m_2,...,\vec m_k,\vec n_1,\vec n_2,...,\vec n_k\in\mbb{N}^{\:\ell-\ell'}.\eqno(5.25)$$
Set
$$\mbb{B}_l^{(i)}=\mbb{A}\td{\ptl}^{\vec m_i},\;\;\;\mbb{B}_r^{(i)}=\td{\ptl}^{\vec n_i}\mbb{A}\qquad\for\;\;i\in\ol{1,k}.\eqno(5.26)$$
In this case
$$\mbb{B}_r^{(i)}\mbb{B}_l^{(j)}=\td{\ptl}^{\vec n_i}\mbb{A}\td{\ptl}^{\vec m_j}\qquad\for\;\;i,j\in\ol{1,k}.\eqno(5.27)$$
Set
$$|\vec m|=\sum_{i=\ell'+1}^\ell m_i\qquad\for\;\;\vec m=(m_{\ell'+1},...,m_{\ell})\in\mbb{N}^{\:\ell-\ell'}.\eqno(5.28)$$
Moreover, we assume
$$|\vec m_1|+|\vec n_1|\equiv |\vec m_2|+|\vec n_2|\equiv\cdots\equiv|\vec m_k|+|\vec n_k|\;\;(\mbox{mod}\:2).\eqno(5.29)$$
We define the linear map $\rho:\hat\mbb{A}\rta\hat\mbb{A}$ by
$$\rho(\td{\ptl}^{\vec n_i}a\ptl^{\vec m_j})=(-1)^{|\vec n_i|+|\vec m_j|}\td{\ptl}^{\vec n_j}\tau(a)\ptl^{\vec m_i}\qquad\for\;\;i,j\in\ol{1,k},\;\td{\ptl}^{\vec n_i}a\ptl^{\vec m_j}\in \mbb{B}_r^{(i)}\mbb{B}_l^{(j)}.\eqno(5.30)$$
\pse

Set
$$\Psi=\left\{\left(\begin{array}{cccc}a_{1,1}& a_{1,2}&\cdots& a_{1,k}\\ a_{2,1}&a_{2,2}&\cdots& a_{2,k}\\ \vdots&\vdots& &\vdots\\  a_{k,1}&a_{k,2}&\cdots& a_{k,k}\end{array}\right)\mid a_{i,j}\in\mbb{B}_r^{(i)}\mbb{B}_l^{(j)};\;i,j\in\ol{1,k}\right\}.\eqno(5.31)$$
Then $\Psi$ is an associative subalgebra of $M_{k\times k}(\mbb{A})$. Moreover, we define a linear map $\ast:\Psi\rta\Psi$ by
$$\left(\begin{array}{cccc}a_{1,1}& a_{1,2}&\cdots& a_{1,k}\\ a_{2,1}&a_{2,2}&\cdots& a_{2,k}\\ \vdots&\vdots& &\vdots\\  a_{k,1}&a_{k,2}&\cdots& a_{k,k}\end{array}\right)^{\ast}=\left(\begin{array}{cccc}\rho(a_{1,1})& \rho(a_{2,1})&\cdots& \rho(a_{k,1})\\ \rho(a_{1,2})&\rho(a_{2,2})&\cdots& \rho(a_{k,2})\\ \vdots&\vdots& &\vdots\\ \rho(a_{1,k})&\rho(a_{2,k})&\cdots& \rho(a_{k,k})\end{array}\right).\eqno(5.32)$$
It can be verified that $\ast$ is an involutive anti-automorphism of $\Psi$.
Define
$$o(\vec\mbb{B}_l,\vec\mbb{B}_r)=\{A\in\Psi\mid A^{\ast}=-A\}.\eqno(5.33)$$
Then $o(\vec\mbb{B}_l,\vec\mbb{B}_r)$ forms a Lie subalgebra of $M_{k\times k}(\mbb{A})$ with the Lie bracket defined in (5.1). 
\psp

{\bf Theorem 5.2 (Benkart, Xu and Zhao, [BXZ])}. {\it The Lie algebra} $o(\vec\mbb{B}_l,\vec\mbb{B}_r)$ {\it is simple}.
\psp

Next we set
$$\Psi'=\left\{\left(\begin{array}{ccc}a_{1,1}&\cdots& a_{1,2k}\\ \vdots& &\vdots\\  a_{2k,1}&\cdots& a_{2k,2k}\end{array}\right)\mid a_{i,j},a_{k+i,j},a_{i,k+j},a_{k+i,k+j}\in\mbb{B}_r^{(i)}\mbb{B}_l^{(j)};\;i,j\in\ol{1,k}\right\}.\eqno(5.34)$$
Then $\Psi'$ is an associative subalgebra of $M_{2k\times 2k}(\mbb{A})$. Moreover, we define a linear map $\ast:\Psi'\rta\Psi'$ by
$$\left(\begin{array}{cccc}a_{1,1}& a_{1,2}&\cdots& a_{1,2k}\\ a_{2,1}&a_{2,2}&\cdots& a_{2,2k}\\ \vdots&\vdots& &\vdots\\  a_{2k,1}&a_{2k,2}&\cdots& a_{2k,2k}\end{array}\right)^{\ast}=\left(\begin{array}{cccc}\rho(a_{1,1})& \rho(a_{2,1})&\cdots& \rho(a_{2k,1})\\ \rho(a_{1,2})&\rho(a_{2,2})&\cdots& \rho(a_{2k,2})\\ \vdots&\vdots& &\vdots\\ \rho(a_{1,2k})&\rho(a_{2,2k})&\cdots& \rho(a_{2k,2k})\end{array}\right).\eqno(5.35)$$
It can be verified that $\ast$ is an involutive anti-automorphism of $\Psi'$.

Denote
$$S=\left(\begin{array}{cc}&I_k\\ -I_k&\end{array}\right),\eqno(5.36)$$
where $I_k$ is the $k\times k$ identity matrix. Moreover,
we define
$$\rho(A)=-SA^{\ast}S\qquad\for\;\;A\in \Psi'.\eqno(5.37)$$
Then $\rho$ is involutive anti-automorphism of $\Psi'$. We set
$$sp(2k,\vec\mbb{B}_l,\vec\mbb{B}_r)=\{A\in\Psi'\mid \rho(A)=-A\}.\eqno(5.38)$$
The subspace $sp(2k,\vec\mbb{B}_l,\vec\mbb{B}_r)$ forms a Lie subalgebra of $M_{2k\times 2k}(\mbb{A})$ with the Lie bracket defined in (5.1).
\psp

{\bf Theorem 5.3 (Benkart, Xu and Zhao, [BXZ])}. {\it The Lie algebra} $sp(2k,\vec\mbb{B}_l,\vec\mbb{B}_r)$ {\it is simple}.
\pse

We remark that in [BXZ], we obtained Theorems 5.1-5.3 for any simple associative algebra satisfying certain properties as $\mbb{A}$ do. These simple algebras  are not finitely-graded. All the classical Lie superalgebras over generalized Weyl algebras had been obtained.

\vspace{1cm}

\noindent{\Large \bf References}

\hspace{0.5cm}

\begin{description}

\item[{[BXZ]}] G. Benkart, X. Xu and K. Zhao, Classical Lie superalgebras over simple associative algebras, {\it Preprint}.

\item[{[Bl]}] R. Block, On torsion-free abelian groups and Lie algebras,
{\it Proc.
Amer. Math. Soc.} {\bf 9} (1958), 613-620.

\item[{[BD]}] A. Beilinson and V. Drinfeld, {\it Chiral Algebras I}, preprint.

\item[{[BDK]}] B. Bakalov, A. D'Andrea and V. Kac, Theory of finite pseudoalgebras,  ArXiv:Math-QA/0007121.

\item[{[Bo]}] R. E. Borcherds, Vertex algebras, Kac-Moody algebras, and the
Monster,
{\it Proc. Natl. Acad. Sci. USA} {\bf 83} (1986), 3068-3071.

\item[{[BP]}] J. Bergen and D. S. Passman, Simple Lie algebras of Special type,
{\it J. Algebra} {\bf 227} (2000), 45-57.

\item[{[CK]}] S.-J. Cheng and V. Kac, Structure of some $\mbb{Z}$-graded Lie
superalgebras
of vector fields, {\it Transformation Groups} {\bf 4} (1999).

\item[{[DZ1]}]  D. Djokovic and K. Zhao, Derivations, isomorphisms and second cohomology of generalized Block algebras, {\it Algebra Colloq.} {\bf 3} (1996), 245-272.

\item[{[DZ2]}] ---, Some simple subalgebras of generalized
Block algebras, {\it J. Algebra} {\bf 192} (1997), 74-101.

\item[{[DZ3]}]  ---, Generalized Cartan type $W$ Lie algebras in characteristic 0, {\it J. Algebra} {\bf 195} (1997), 170-210.

\item[{[DZ4]}] ---, Generalized Cartan type S Lie algebras in characteristic
zero,
{\it J. Algebra} {\bf 193} (1997), 144-179.

\item[{[DZ5]}] ---, Some infinite-dimensional simple Lie algebras in
characteristic
0 related to those of Block, {\it J. Pure Appl. Algebra} {\bf 127} (1998),
153-165.

\item[{[FLM]}] I. Frenkel, J. Lepowsky and A. Meurman, {\it Vertex Operator
Algebras and the Monster}, Pure and Applied Math. Academic Press, 1988.

\item[{[J1]}] D. A. Jordan, Simple Lie rings of derivations of commutative rings, {\it J. London Math. Soc. (2)} {\bf 18} (1978), 443-448. 

\item[{[J2]}] ---, On ideals of a Lie algebra of derivations, {\it J. London Math. Soc. (2)} {\bf 33} (1986), 33-39.

\item[{[J3]}] ---, On the simplicity of Lie algebras of derivations of commutative algebras, {\it J. Algebra} {\bf 228} (2000), 580-585.

\item[{[K1]}] V. G. Kac, Simple graded Lie algebras of finite growth,  {\it
Funct.
Anal. Appl.} {\bf 1} (1967), 328-329.

\item[{[K2]}] ---, A description of filtered Lie algebras whose associated
graded Lie
algebras are of Cartan types, {\it Math. of USSR-Izvestijia} {\bf 8} (1974),
801-835.

\item[{[K3]}] ---, Lie superalgebras, {\it Adv. Math.} {\bf 26} (1977),
8-96.

\item[{[K4]}] ---, {\it Infinite dimensional Lie algebras}, Third Edition,
Cambridge
University Press, 1990.

\item[{[K5]}] ---, {\it Vertex algebras for beginners}, University
lectures
series, Vol {\bf 10}, AMS. Providence RI, 1996.

\item[{[K6]}] ---, Classification of infinite-dimensional simple linearly
compact Lie
superalgebras, {\it Adv. Math.} {\bf 139} (1998), 1-55.

\item[{[Ki]}] I. L. Kantor, Jordan and Lie superalgebras determined by a Poisson algebra, {\it Algebra and Analysis (Tomsk, 1989)}, 55-80, Amer. Math. Soc. Ser. 2, 151, {\it Amer. Math. Soc., Providence, RI, 1992}.

\item[{[Kn]}] N. Kawamoto, Generalizations of Witt algebras over a field of
characteristic
zero,
{\it Hiroshima Math. J.} {\bf 16} (1986), 417-462.

\item[{[L1]}] D. A. Leites, New Lie superalgebras and mechanics, (Russian) {\it Dokl. Akad. Nauk SSSR} {\bf 236} (1977), 804-807.

\item[{[L2]}] ---, Automorphisms and real forms of simple Lie superalgebras of formal vector fields, (Russian) {\it Problems in group theory and homological algebra}, 126-128, 139, {\it Yaroslav. Gos. Univ. Yaroslav,} 1983.

\item[{[Mo1]}] O. Mathieu, Classification of simple graded Lie algebras of growth $\leq 1$, {\it Invent. Math.} {\bf 86} (1986), 371-426.

\item[{[Mo2]}] ---, Classification of simple graded Lie algebras of finite growth, {\it Invent. Math.} {\bf 108} (1992), 455-519.

\item[{[Mr]}] R. V. Moody, Lie algebras associated with generalized Cartan
matrices,
{\it Bull. Amer. Math. Soc.} {\bf 73} (1967), 217-221.

\item[{[O]}] J. Marshall Osborn, New simple infinite-dimensional Lie
algebras of
characteristic 0, {\it J. Algebra} {\bf 185} (1996), 820-835.

\item[{[OZ1]}] J. Marshall Osborn and K. Zhao, Generalized Cartan type K Lie
algebras in characteristic 0, {\it Commun. Algebra} {\bf 25} (1997),
3325-3360.

\item[{[OZ2]}] ---, Generalized Poisson brackets and Lie algebras of type H
in
characteristic 0, {\it Math. Z.} {\bf 230} (1999), 107-143.

\item[{[OZ3]}] ---, Infinite dimensional Lie algebras of generalized Block
type,
{\it Proc. Amer. Math. Soc.} {\bf 127} (1999), 1641-1650.

\item[{[OZ4]}] ---, Doubly $\mbb{Z}$-graded Lie algebras containing a
Virasoro algebra,
{\it J. Algebra} {\bf 219} (1999), 266-298.

\item[{[P1]}] D. S. Passman, Simple Lie algebras of Witt type, {\it J.
Algebra} {\bf 206}
(1998), 682-692.

\item[{[P2]}] ---, Simple Lie color algebras of Witt type, {\it J. Algebra}
{\bf 208}
(1998), 698-721.

\item[{[S1]}] I. Shchepochkina, Exceptional simple infinite-dimensional Lie
superalgebras,
{\it C. R. Acad. Bulgare Sci.} {\bf 36} (1983), 313-314.

\item[{[S2]}] ---, Five simple exceptional Lie superalgebras of vector
fields, {\it Funt. Anal. Appl.} {\bf 33} (1999), 208-219.

\item[{[S3]}] ---, The five exceptional simple Lie superalgebras of vector
fields and their fourteen regradings, {\it Represent. Theory} {\bf 3} (1999), 373-415.

\item[{[SX1]}] Y. Su and X. Xu, Structure of divergence-free Lie algebras,  ArXiv:Math.QA/0005286.

\item[{[SX2]}] ---, Central simple Poisson algebras,  ArXiv:Math.QA/0011086.

\item[{[SXZ]}] Y. Su, X. Xu and H. Zhang, Derivation-Simple algebras and the
structures of generalized Lie algebras of Witt type,  {\it J. Algebra} {\bf 233} (2000), 642-662. 

\item[{[W]}] B. Weisfeiler, Infinite dimensional filtered Lie algebras and their relation with graded Lie algebras, {\it Funct. Anal. i Prilozhen} {\bf 2} (1968), 94-95.

\item[{[X1]}] X. Xu, On simple Novikov algebras and their irreducible modules, {\it J. Algebra} {\bf 185} (1996), 905-934.

\item[{[X2]}] ---, Novikov-Poisson algebras, {\it J. Algebra} {\bf 190} (1997), 253-279.

\item[{[X3]}] ---, {\it Introduction to Vertex Operator Superalgebras and Their Modules},
Kluwer Academic Publishers, Dordrecht/Boston/London, 1998.

\item[{[X4]}] ---, Generalizations of Block algebras, {\it Manuscripta
Math.} {\bf 100} (1999), 489-518.

\item[{[X5]}] ---, Variational calculus of supervariables and related algebraic structures, {\it J. Algebra} {\bf 223} (2000), 396-437.

\item[{[X6]}] ---,  Quadratic conformal superalgebras, {\it J. Algebra} {\bf 231} (2000), 1-38.

\item[{[X7]}] ---, New generalized simple Lie algebras of Cartan type over
a field with
characteristic 0, {\it J. Algebra} {\bf 224} (2000), 23-58.

\item[{[X8]}] ---, Simple conformal superalgebras to finite growth, {\it Algebra Colloquium} {\bf 7} (2000), 205-240.

\item[{[X9]}] ---, Equivalence of conformal superalgebras to Hamiltonian superoperators, to appear in {\it Algebra Colloquium} {\bf 8} (2001), 63-92.

\item[{[Z1]}] K. Zhao, Generalized Cartan type S Lie algebras in
characteristic zero (II),
{\it Pacific J. Math.} {\bf 192} (2000), 431-454.

\item[{[Z2]}] ---, A class of infinite-dimensional simple Lie algebras, {\it J. London Math. Soc. (2)}  {\bf 62} (2000), 71-84.

\end{description}
\end{document}